\documentstyle[fleqn,12pt]{article}
\begin{document}
\pagenumbering{arabic} \setcounter{page}{1} \pagestyle{plain}
\baselineskip=16pt

\thispagestyle{empty}
\begin{center} {\Large\bf The Hopf algebra
structure of the Z$_3$-graded quantum supergroup GL$_{q,j}(1\vert
1)$ }
\end{center}

\vspace{1cm} \noindent Salih Celik\footnote{Electronic mail:
sacelik@yildiz.edu.tr}, Erg\"un Yasar\footnote{Electronic mail:
ergyasmat@gmail.com }

\noindent Yildiz Technical University, Department of Mathematics,
34210 Davutpasa-Esenler, Istanbul, Turkey.

\vspace{2cm}

In this work, we give some features of the Z$_3$-graded quantum
supergroup.

\vspace{2cm}\noindent {\bf I. \, INTRODUCTION}

Recently, there have been many attempts to generalize Z$_2$-graded
constructions to the Z$_3$-graded case$^{1-3}$. The Z$_3$-graded
quantum space that generalizes the Z$_2$-graded space called a
superspace$^4$, was studied using the methods of Ref. 5. The first
author studied the noncommutative geometry of the Z$_3$-graded
quantum superplane$^6$. In this work, we have investigated the
Hopf algebra structure of the Z$_3$-graded quantum supergroup
GL$_{q,j}(1\vert 1)$.

Let us shortly investigate a general Z$_3$-graded algebraic
structure. Let $z$ be a Z$_3$-graded  variable. Then we say that
the variable $z$ satisfies the relation
$$z^3 = 0.$$
If $f(z)$ is an arbitrary function of the variable $z$, then the
function $f(z)$ becomes a polynomial of degree two in $z$, that
is,
$$f(z) = a_0 + a_1 z + a_2 z^2,$$
where $a_0$, $a_2$, $a_1$ denote three fixed numbers whose grades
are $grad (a_0) = 0$, $grad (a_2) = 1$ and $grad(a_1) = 2$,
respectively.

The cyclic  group Z$_3$ can be represented in the complex plane by
means of the cubic roots of 1: let $j = e^{{2\pi i}\over 3}$ $(i^2
= - 1)$. Then one has
$$j^3 = 1 \quad \mbox{and} \quad j^2 + j +1 = 0.  $$
One can define the Z$_3$-graded commutator $[A,B]$ as
$$[A,B]_{Z_3} = AB - j^{ab} BA, $$
where $grad(A) = a$ and $grad(B) = b$. If $A$ and $B$ are
$j$-commutative, then we have
$$AB = j^{ab} BA.$$

\noindent {\bf II. \,  REVIEW OF THE ALGEBRA OF FUNCTIONS ON THE
Z$_3$-GRADED QUANTUM SUPERPLANE}

The Z$_3$-graded quantum superplane is defined as an associative
unital algebra generate by $x$ and $\theta$ satisfying$^6$
$$ x \theta = q \theta x, \qquad \theta^3 = 0 \eqno(1)$$
where $q$ is a nonzero complex deformation parameter. Here, the
coordinate $x$ with respect to the Z$_3$-grading is of grade 0 and
the coordinate $\theta$ with respect to the Z$_3$-grading is of
grade 1. This associative algebra over the complex numbers is
known as the algebra of polynomials over the quantum superplane
and we shall denote it by ${\cal R}_q (1\vert 1)$, that is,
$$ {\cal R}_q (1\vert 1) \ni \left(\matrix{x \cr \theta}\right)
\qquad \Longleftrightarrow \qquad x \theta = q \theta x, \qquad
\theta^3 = 0.$$ If we denote the dual of the ${\cal R}_q (1\vert
1)$ by ${\cal R^\star}_{q,j} (1\vert 1)$, one has
$$ {\cal R^\star}_{q,j} (1\vert 1) \ni \left(\matrix{\varphi
\cr \ y}\right)  \quad \Longleftrightarrow \qquad \varphi\ y = q j
y \varphi, \quad \varphi^3 = 0.\eqno(2)$$ Here,
$$ [{\cal R}_q (1\vert 1)]^\star = {\cal R^\star}_{q,j} (1\vert
1).$$ We define the extended quantum superplane to be the algebra
that contains ${\cal R}_q (1\vert 1)$, the unit and $x^{-1}$, the
inverse of $x$, which obeys
$$x x^{-1} = 1 =  x^{-1} x.$$
We denote the extended algebra by ${\cal A}$. We know that the
algebra ${\cal A}$ is a Z$_3$-graded Hopf algebra.$^6$

\newpage\noindent {\bf III. \, A PERSPECTIVE TO Z$_3$-GRADED
$h$-DEFORMATION }

We know that the commutation relation between the coordinate $x'$
and the coordinate $\theta'$ of the Z$_3$-graded quantum
superplane is in the form
$$ x' \theta' - q \theta' x' = 0. $$
We now introduce new coordinates $x$ and $\theta$, in terms of
$x'$ and $\theta'$ as
$$ x = x',\qquad \theta = \theta' - \frac{h}{q-1} x' \eqno(3)$$
as in Ref. 7. This transformation is singular in the $q
\rightarrow 1$ limit. Using relation (1), it is easy to verify
that
$$ x \theta = q \theta x + h x^2 \eqno(4)$$
where the new deformation parameter $h$ commutes with the
coordinate $x$. Also, since the grassmann coordinate $\theta'$
satisfies
$$\theta'^3 = 0$$
one obtains
$$\theta^3 = 0 \eqno(5)$$
provided that
$$\theta h = q j h \theta, \qquad h^3 = 0. \eqno(6)$$
Taking the $q \rightarrow 1$ limit we obtain the following
relations which define the Z$_3$-graded $h$-superplane
$$ x \theta =  \theta x + h x^2, \qquad \theta^3 = 0.\eqno(7)$$
Also, it can be obtained the Z$_3$-graded $h$-supergroup with
Aghamohammadi's approach in Ref. 7. So, it can be investigated the
differential geometry of this group. This work is in progress.

\noindent {\bf IV. \, Z$_3$-GRADED QUANTUM SUPERGROUPS }

\noindent {\bf A. \, Quantum matrices in Z$_3$-graded superspace}

In this section, we shall consider the Z$_3$-graded structures of
the quantum 2x2 supermatrices. We know, from section 2, that the
Z$_3$-graded quantum superplane ${\cal R}_q (1\vert 1)$ is
generated by coordinates $x$ and $\theta$, with the commutation
rules (1), and the dual Z$_3$-graded quantum superplane ${\cal
R^\star}_{q,j} (1\vert 1)$ as generated by $\varphi$ and $y$ with
the relations (2).

Let $T$ be a 2x2 (super)matrix in Z$_3$-graded space,
$$T = \left(\matrix{a & \beta \cr \gamma & d}\right) \eqno(8)$$
where $a$ and $d$ with respect to the Z$_3$-grading are of grade
0, and $\beta$ and $\gamma$ with respect to the Z$_3$-grading are
of grade 2 and of grade 1, respectively. We now consider linear
transformations with the following properties:
$$\ T : {\cal R}_q
(1\vert 1) \longrightarrow {\cal R}_q (1\vert 1),\qquad
  T : {\cal R^\star}_{q,j} (1\vert 1) \longrightarrow
{\cal R^\star}_{q,j} (1\vert 1). \eqno(9)$$ We assume that the
entries of $T$ are $j$-commutative with the elements of ${\cal
R}_q (1\vert 1)$ and ${\cal R^\star}_{q,j}(1\vert 1)$, i.e. for
example,
$$a x = x a, \qquad \theta \beta = j^2 \beta \theta, $$
etc. As a consequence of the linear transformations in (9) the
elements
$$\tilde{x} = a x + \beta \theta, \qquad \tilde{\theta} = \gamma x + d \theta \eqno(10)$$
should satisfy the relations (1):
$$ \tilde{x} \tilde{\theta} = q \tilde{\theta} \tilde{x}, \qquad \tilde{\theta}^3 = 0.$$
Using these relations , one has
$$a \gamma = q \gamma a, \qquad d \gamma = q  \gamma d,$$
$$d \beta = q^{-1} j \beta d,\qquad \gamma^3 = 0.$$
Similarly, the elements
$$\tilde{\varphi} = a \varphi + j^2\beta y, \qquad \tilde{y} = j \gamma \varphi + d y \eqno(11)$$
must be satisfy the relations (2). Using these relations , one has
$$a \beta = q^{-1} j^{-1} \beta a, \qquad \beta^{3} = 0. $$ Also, if we
use the following relation in Ref. 6 (see, page 4262, eq. (19))
$$ \tilde{x} \tilde{y} = q \tilde{y} \tilde{x} + (j^{2}-1) \tilde{\varphi} \tilde{\theta},$$
we have
$$a d = d a + q^{-1} (1 - j) \beta \gamma, \qquad
  \beta \gamma = q^2 \gamma \beta.$$
Consequently, we have the following commutation relations between
the matrix elements of $T$ which is given in Ref. 6:
$$a \beta = q^{-1} j^{-1} \beta a, \qquad d \beta = q^{-1} j \beta d, $$
$$a \gamma = q \gamma a, \qquad d \gamma = q  \gamma d,$$
$$a d = d a + q^{-1} (1 - j) \beta \gamma, \qquad
  \beta \gamma = q^2 \gamma \beta,$$
$$ \beta^{3} = 0, \qquad \gamma^3 = 0. \eqno(12)$$
The Z$_3$-graded quantum (super)determinant is defined by
$$ D_{q,j}(T)= a d^{-1} + a d^{-1}\gamma a^{-1} \beta d^{-1}
+ a d^{-1}\gamma a^{-1} \beta d^{-1}\gamma a^{-1} \beta d^{-1}
\eqno(13)$$ provided $d$ and $a$ are invertible. The commutation
relations between the matrix elements of $T$ and its
(super)determinant:
$$ a D_{q,j}(T) = D_{q,j}(T) a, \quad \beta D_{q,j}(T) = j^{2} D_{q,j}(T) \beta  $$
$$ \gamma D_{q,j}(T) = D_{q,j}(T) \gamma, \quad d a D_{q,j}(T) = D_{q,j}(T) a d  \eqno(14) $$
If we take
$$\Delta_{1}= a d - q^{-1} \beta \gamma \qquad \Delta_{2}= d a - q j \gamma
\beta$$ the Z$_3$-graded quantum (super)inverse of $T$ becomes
$$T^{-1} = \left(\matrix{d {\Delta_{1}}^{-1}& - q j \beta {\Delta_{2}}^{-1}
\cr - q^{-1}\gamma {\Delta_{1}}^{-1} & a
{\Delta_{2}}^{-1}}\right).$$ Using the relations (12), one can
check that the following relations:
$$ \Delta_{1} a = a \Delta_{1}, \quad \Delta_{1} d = d \Delta_{1}$$
$$ \Delta_{2} a = a \Delta_{2}, \quad \Delta_{2} d = d \Delta_{2}$$
$$ \Delta_{k} \beta = q^{-2} \beta \Delta_{k}, \quad \Delta_{k} \gamma = q^{2} \gamma \Delta_{k}, \quad k=1,2$$
$$\Delta_{1} \Delta_{2} = \Delta_{2} \Delta_{1}. \eqno(15)$$
The Z$_3$-graded quantum (super)determinant of $T$, according to
these facts, is given by
$$ D_{q,j}(T)= a^{2} {\Delta_{2}}^{-1}.$$
Of course, the Z$_3$-graded quantum (super)determinant of $T^{-1}$
may also be defined and it is of the form
$$ D_{q,j}(T^{-1})= d^{2} {\Delta_{1}}^{-1}.$$ Explicitly,
$$ D_{q,j}(T^{-1})= d a^{-1} + d a^{-1}\beta d^{-1} \gamma a^{-1}
+ d a^{-1}\beta d^{-1} \gamma a^{-1}\beta d^{-1} \gamma
a^{-1}.\eqno(16)$$ Let's now consider the multiplication of two
Z$_3$-graded quantum (super) matrices. If we take them as
$$T_{1} = \left(\matrix{a_{1} & \beta_{1} \cr \gamma_{1} & d_{1}}\right)
\quad \mbox{and} \quad T_{2} = \left(\matrix{a_{2} & \beta_{2} \cr
\gamma_{2} & d_{2}}\right)$$ where the matrix elements of $T_{1}$
and $T_{2}$ satisfy the relations (12). Then the matrix elements
of $T_{1}T_{2}$ leave invariant the relations (12), as expected.
Here, we assume that the commutation relations between the
elements of $T_{1}$ and $T_{2}$ are as follows
$$\beta_{1} \gamma_{2}= j \gamma_{2}\beta_{1}, \quad \beta_{1} \beta_{2}= j^{2} \beta_{2}\beta_{1},$$
$$\gamma_{1} \beta_{2}= j \beta_{2}\gamma_{1}, \quad \gamma_{1} \gamma_{2}= j^{2} \gamma_{2}\gamma_{1}$$
and the elements whose gradings are 0 commute with all the other
elements. Also, the Z$_3$-graded quantum (super)determinant is not
central of the Z$_3$-graded quantum (super)group, although the
Z$_2$-graded quantum superdeterminant is central of the
Z$_2$-graded quantum supergroup with two parameters.$^8$

We shall denote with GL$_{q,j}(1\vert 1)$ the quantum supergroup
in Z$_3$-graded space determined by generators $a$, $\beta$,
$\gamma$, $d$ satisfying the commutation relations (12).

Also, we can define the Z$_3$-graded quantum (super)transpoze of
$T$ as
$$T^{st} = \left(\matrix{a & \gamma \cr j \beta & d}\right). \eqno(17)$$
Note that, the transformations
$$\ T^{st}: {\cal R}_q
(1\vert 1) \longrightarrow {\cal R}_q (1\vert 1),\qquad
  T^{st} : {\cal R^\star}_{q,j} (1\vert 1) \longrightarrow
{\cal R^\star}_{q,j} (1\vert 1) \eqno(18)$$ with together the
transformations (9), will give the relations (12). Here, the
action of $T^{st}$ on the coordinate functions as follows
$$ \hat{x} = x a + j \theta \beta, \quad \hat{\theta} = x \gamma + \theta d. $$
One can show that $T^{-1}\in GL_{q^{-1},j^{-1}}(1\vert 1)$.
Indeed, if we denote the Z$_3$-graded quantum (super)inverse of
$T$ with
$$T^{-1} = \left(\matrix{A & B \cr C & D}\right) \eqno(19) $$
then we have
$$A = a^{-1} + a^{-1} \beta d^{-1}\gamma a^{-1} + a^{-1} \beta d^{-1}\gamma a^{-1} \beta d^{-1}\gamma a^{-1}, $$
$$B = -a^{-1} \beta d^{-1} -a^{-1} \beta d^{-1} \gamma a^{-1} \beta
d^{-1} -a^{-1} \beta d^{-1} \gamma a^{-1} \beta d^{-1} \gamma
a^{-1} \beta d^{-1},$$
$$C = -d^{-1} \gamma a^{-1} -d^{-1} \gamma
a^{-1}\beta d^{-1} \gamma a^{-1},$$
$$D = d^{-1} + d^{-1} \gamma a^{-1}
\beta d^{-1} + d^{-1} \gamma a^{-1} \beta d^{-1}\gamma a^{-1}
\beta d^{-1}.\eqno(20) $$ Now, using the relations (15) and (12)
one has $T^{-1}\in GL_{q^{-1},j^{-1}}(1\vert 1)$.

The relations between the matrix elements of $T$ and
(super)inverse of $T$ are important, because we will use them to
set up the Z$_3$-graded differential geometric structure of the
Z$_3$-graded quantum supergroup.$^9$ These relations are as
follows:
$$ a A = j^{2} A a + 1-j^{2}, \quad a B = q^{-1} j^{2} B a, $$
$$ a C = q C a, \quad a D = D a, $$
$$ \beta A = q^{-1} j^{2} A \beta, \quad \beta B = q ^{-2} B \beta, $$
$$ \beta C = C \beta, \quad \beta D = q^{-1} j D \beta, $$
$$ \gamma A = q A \gamma, \quad \gamma B = B \gamma, $$
$$ \gamma C = q^{2} C \gamma, \quad \gamma D = q D \gamma, $$
$$ d A = A d, \quad d B = q^{-1} j B d, $$
$$ d C = q C d, \quad d D = j D d + 1-j.\eqno(21) $$

Also, it can be investigated some properties of the quantum
(super)matrices in the quantum supergroup GL$_{q,j}(1\vert 1)$.
So, perhaps any element of $GL_{q,j}(1\vert 1)$ can be expressed
as the exponential of a matrix of non-commuting elements, like the
group $GL_q(1\vert1)$. This work is also in progress.

\noindent {\bf B. \, Hopf algebra structure of the Z$_3$-graded
GL$_{q,j}(1\vert 1)$ }

In this section, we shall build up the Hopf algebra structure of
the Z$_3$-graded quantum supergroup GL$_{q,j}(1\vert 1)$. For
this, we shall introduce three operators $\Delta$, $\epsilon$ and
$S$ on the GL$_{q,j}(1\vert 1)$, which are called the coproduct
(comultiplication), the counit and the coinverse (antipode),
respectively. The coproduct
$$\Delta : GL_{q,j}(1\vert
1) \longrightarrow \mbox{GL}_{q,j}(1\vert 1) \otimes
GL_{q,j}(1\vert 1)$$ is defined by
$$\Delta ( T ) = T \dot{\otimes} T \eqno(22)$$
where $\dot{\otimes}$ stands for the usual tensor product and the
dot refers to the summation over repeated indices and reminds us
about the usual matrix multiplication. The coproduct $\Delta$ is
an algebra homomorphism which is co-associative, that is
$$(\Delta \otimes { \mbox{id}}) \circ \Delta = ({ \mbox {id}} \otimes \Delta) \circ \Delta \eqno(23) $$
where $\circ$ stands for the composition of maps and id denotes
the identity mapping. Also, the multiplication in the algebra of
matrix entries of $T$ is defined as
$$(A \otimes B) (C \otimes D) = j^{grad(B) grad(C)} AC \otimes BD.$$

\noindent The action on the generators of the GL$_{q,j}(1\vert 1)$
of $\Delta$ is
$$\Delta( a ) = a \otimes a + \beta \otimes \gamma, \quad
\Delta( \beta ) = a \otimes \beta + \beta \otimes d, $$
$$\Delta( \gamma ) = \gamma \otimes a + d \otimes \gamma, \quad
\Delta( d ) = \gamma \otimes \beta + d \otimes d \eqno(24) $$
where $\otimes$ denotes the tensor product. Of course, the
coproduct $\Delta$ leaves invariant the relations (12). The counit
$$\epsilon: GL_{q,j}(1\vert
1) \longrightarrow {\cal C}$$ is defined by
$$ \epsilon ( T ) = I . \eqno(25)$$  The action on the generators of
$GL_{q,j}(1\vert 1)$ of $\epsilon$ is
$$\epsilon(a) = 1, \qquad \epsilon(\beta) = 0, \qquad \epsilon(\gamma) = 0, \qquad \epsilon(d) = 1.
\eqno(26)$$ The counit $\epsilon$ is an algebra homomorphism such
that
$$\mu \circ (\epsilon \otimes {\mbox{id}}) \circ \Delta
  = \mu' \circ ({\mbox{id}} \otimes \epsilon) \circ \Delta \eqno(27)$$
where
$$\mu : {\cal C} \otimes {GL_{q,j}(1\vert 1)} \longrightarrow {GL_{q,j}(1\vert 1)}, \qquad
  \mu' : {GL_{q,j}(1\vert 1)} \otimes {\cal C} \longrightarrow {GL_{q,j}(1\vert 1)} $$
are the canonical isomorphisms, defined by
$$\mu(c \otimes a) = c a = \mu'(a \otimes c), \qquad \forall a \in {GL_{q,j}(1\vert 1)},
  \quad \forall c \in {\cal C}. $$
Thus, we have verified that GL$_{q,j}(1\vert 1)$ is a bialgebra
with the multiplication $m$ satisfying the associativitiy axiom
$$ m \circ ( m \otimes \mbox{id}) = m \circ ( \mbox{id} \otimes m) $$
where $m(a \otimes b) = ab.$

\noindent A bialgebra with the extra structure of the coinverse is
called a Hopf algebra.$^{10}$

\noindent The coinverse
$$ S: GL_{q,j}(1\vert
1) \longrightarrow GL_{q,j}(1\vert 1)$$ is defined by
$$ S (T) =
T^{-1}.\eqno(28)$$ The coinverse $S$ is an algebra
anti-homomorphism which satisfies
$$m \circ (S \otimes {\mbox{id}}) \circ \Delta = \epsilon
  = m \circ ({\mbox{id}} \otimes S) \circ \Delta.\eqno(29) $$
The coproduct, counit and coinverse which are specified above
supply $GL_{q,j}(1\vert 1)$ with a Hopf algebra structure. It can
be show that $GL_{q,j}(1\vert 1)$ has a Z$_3$-graded differential
geometric structure.$^9$

\noindent {\bf ACKNOWLEDGMENT}

This work was supported in part by TBTAK the Turkish Scientific
and Technical Research Council.

\baselineskip=12pt

{\small



$^1$ R. Kerner and V. Abramov, Rep. Math. Phys. {\bf 43}, 179
(1999).

$^2$ B. Le Roy, J. Math. Phys. {\bf 37}, 474 (1996).

$^3$ V. Abramov and N. Bazunova, math-ph/0001041, (2001).

$^4$ W. S. Chung, J. Math. Phys. {\bf 35}, 2497 (1993).

$^5$ J. Wess and B. Zumino, Nucl. Phys. B (Proc. Suppl.) {\bf 18}
B, 302 (1990).

$^6$ S. Celik, J.Phys. A: Math. Gen. {\bf 35}, 4257 (2002).

$^7$ A. Aghamohammadi, M. Khorrami and A. Shariati, A. J. Phys. A:
Math. Gen. {\bf 28}, 225 (1995);

L. Dabrowski and P. Parashar, Lett. Math. Phys. {\bf 38}, 331
(1996).

$^8$ L. Dabrowski and L. Wang: Phys. Lett. B {\bf 266}, 51 (1991).

$^9$ E. Yasar, "Differential geometry of the Z$_3$-graded quantum
supergroup ", Phd Thesis (in preperation).

$^{10}$ E. Abe, Hopf Algebras, Cambridge Univ. Press, (1980). }

\end{document}